\documentclass[10pt,journal]{IEEEtran}

\usepackage{amsfonts}
\usepackage{amsmath}
\usepackage{amssymb}
\usepackage{amsthm}
\usepackage{hyperref}
\usepackage{graphics}
\usepackage{curves}
\usepackage{epic}
\usepackage{graphicx}
\usepackage{multirow}

\usepackage{floatrow}
\newtheorem{them}{Theorem}

\newtheorem{remk}{Remark}
\newtheorem{exmp}{Example}

\renewcommand{\baselinestretch}{0.98}

\title{A Numerical Approach to Stability of Multi-class Queueing Networks}

\begin{document}

\date{ }

\author{H. Leahu, M. Mandjes (Univ.\ of Amsterdam) \& AM. Oprescu (Vrije Univ.\ Amsterdam)}

\maketitle

\begin{abstract}\noindent
The Multi-class Queueing Network (McQN) arises as a natural multi-class extension of the traditional (single-class) Jackson network.
In a single-class network subcriticality (i.e.\ subunitary nominal workload at every station) entails stability, but this is no longer
sufficient when jobs/customers of different classes (i.e.\ with different service requirements and/or routing scheme) visit the same
server; therefore, analytical conditions for stability of McQNs are lacking, in general.

\noindent In this note we design a numerical (simulation-based) method for determining the stability region of a McQN, in terms of arrival
rate(s). Our method exploits certain (stochastic) monotonicity properties enjoyed by the associated Markovian queue-configuration process.
Stochastic monotonicity is a quite common feature of queueing models and can be easily established in the single-class framework (Jackson
networks); recently, also for a wide class of McQNs, including first-come-first-serve (FCFS) networks, monotonicity properties have been
established. Here, we provide a minimal set of conditions under which the method performs correctly.

\noindent Eventually, we illustrate the use of our numerical method by presenting a set of numerical experiments, covering both single and
multi-class networks.
%For illustrative purposes, for the multi-class examples we test the required monotonicity conditions numerically.
\end{abstract}

\vspace{-2mm}

\section{Introduction}

Multi-class queueing networks (McQNs) provide the mathematical framework for modeling a wide range of stochastic
systems, e.g., manufacturing lines, computer grids and telecommunication systems. They differ from the classical
Jacksonian network model in that the same (physical) item entering the system may require multiple service stages
at the same station, with different service and routing characteristics, thus giving rise to a different class of
jobs. As such, (some) stations behave as multi-class (rather than single-class) queues.

This distinguishing feature has a rather significant impact on the assessment of stability of such networks;
more specifically, while for Jackson networks stability is equivalent to sub-criticality, for some McQNs such
an equivalence does not hold anymore, as demonstrated by a plethora of examples in the literature; see, e.g.
\cite{Bramson:08} for a significant list of examples of subcritical networks which are \emph{not} stable.
It remains true, however, that stability implies subcriticality \cite{Bramson:08}, hence subcriticality is a
necessary, but not sufficient condition for stability.

\newpage

In this note we consider McQNs in which inter-arrival and service times are exponentially distributed; under this
assumption, the queue-configuration process defines a continuous-time Markov chain (Markov process on a discrete
state space) which enables one to employ a more powerful mathematical apparatus. We address the following problem:
\emph{given a certain network, with specified service rates and routing scheme, what is the set of arrival rates
which makes the network stable}? In this context, stability refers to the associated Markovian model, hence positive
Harris recurrence.

While in the Jacksonian framework the answer to the above question is straightforward, under the multi-class paradigm,
in the absence of analytical conditions for stability, one needs to resort to numerical methods. We design a  numerical
(simulation-based) method for solving this problem. Our method, which is among the first schemes of this kind, assumes
some (weak) monotonicity conditions on the associated Markov process, which ensure that the stability region (the set 
of arrival-rate vectors which make the network stable) defines a star-shaped domain in the parameter space. In addition,
the stability region can be recovered by interpolating the boundary points (stability thresholds) in various directions
which, in turn, can be approximated by numerical root-finding methods. Importantly, the required monotonicity conditions 
hold for McQNs in which jobs are executed one at a time; see \cite{LM:2016}.

To test the approach, we performed an extensive set of numerical experiments. We include here a number of illustrative examples.
We show first that the method correctly identifies the predicted stability thresholds when they are available in analytical
form, e.g. for Jackson and Kelly type networks. Furthermore, we apply our numerical method to two instances of multi-class networks
(reentrant lines) where stability conditions are not available, obtaining approximations for the (unknown) stability thresholds.
%New citations 1997 \cite{hasenbein1997necessary} 2005: \cite{gamarnik2005instability} or \cite{zoghby2005modeling}? 
%Finally, we consider an McQN of Kelly type, where thresholds are available in analytic form, and show that our method provides the correct results.

In Section \ref{mcqn:sec} we introduce the mathematical model, the relevant notation and terminology. Furthermore, in Section 
\ref{main:sec} we introduce our method, the necessary assumptions and the (main) convergence result. Finally, in Section
\ref{num:sec} the numerical experiments are presented.

\section{The Mathematical Model}\label{mcqn:sec}

In this section we describe our mathematical model and introduce the notation and terminology which will be used
throughout this note.

\subsection{Multi-class Queueing Networks: The Model}\label{model:sec}

We consider a general McQN model consisting of $\aleph$ stations (each having its own service/queueing policy)
executing $d$ classes of jobs. Each class $k$ is assigned to a specified station $\mathcal{S}(k)$. We further assume
that the mapping $k\longmapsto\mathcal{S}(k)$ is surjective, i.e. each station serves (at least) one class, hence
$1\leq\aleph\leq d$. When the mapping $\mathcal{S}$ is bijective one recovers the standard Jackson Network model.
The set $\{k:\mathcal{S}(k)=i\}$, of all classes assigned to station $i$ will be denoted by $\mathcal{K}_i$.

We now describe the dynamics of the McQN. Jobs of class $k$ enter the network according to a Poisson process with rate
$\theta_k\geq 0$; the case $\theta_k=0$ corresponds to a void arrival process, meaning that class $k$ does not have
external input. Upon arrival, a job of class $k$ is assigned to station $\mathcal{S}(k)$; depending on the underlying
service/queue policy, it either starts receiving service immediately, or it is enqueued in a waiting line. We assume
that jobs of class $k$ require an exponentially distributed service time, with rate $\beta_k>0$, independent of
everything else. After finishing service at station $\mathcal{S}(k)$, a job of class $k$ turns into a job of class
$l$, with probability $R_{kl}$ and moves to station $\mathcal{S}(l)$ (where it follows the corresponding queueing
routine) or leaves the network with probability $R_{k0}:=1-\sum_{l=1}^dR_{kl}$.
To ensure that the network is open, we assume that the matrix $R:=\{R_{kl}\}_{k,l=1,\ldots,d}$ is sub-stochastic,
i.e. $$(I-R)^{-1}=I+R+R^2+\ldots;$$ this condition guarantees that any job will eventually leave the network
(in finite time) with probability one.

An McQN with $R_{k\:k+1}=1$, for $k=1,\ldots,d-1$ and $R_{d0}=1$, such that only class $1$ has non-trivial external
input, i.e., $\theta_2=\ldots=\theta_d=0$, is called a \emph{reentrant line}. Reentrant lines are the most popular
instances of McQNs, as they provide mathematical models for manufacturing systems (assembly lines).

We define the vector of \emph{effective arrival rates} by $$\lambda:=(I-R')^{-1}\theta.$$ Furthermore, the \emph{traffic rate}
(or \emph{nominal workload}) of station $i$ is defined as
\begin{equation}\label{traffic:eq}
\rho_i:=\sum_{k\in\mathcal{K}_i}\frac{\lambda_k}{\beta_k}.
\end{equation}
Station $i$ is called \emph{sub-critical} if $\rho_i<1$ and the network is called sub-critical if every node is.

\newpage

\subsection{The Stability and the Subcriticality Regions}\label{regio:sec}

Under the assumptions in Section \ref{model:sec}, the queue-configuration process defines a Markov process \cite{Dai:95},
$\mathcal{X}:=\{X_t:t\geq 0\}$ (on some suitable state-space $\mathbb{X}$), which depends on the parameter $\theta=(\theta_1,\ldots,\theta_d)\in\Theta$,
where $\Theta\subseteq\mathbb{R}^d_+:=\{\theta\in\mathbb{R}^d:\theta\geq\mathbf{0}\:\text{(componentwise)}\}$ is a pre-specified set;
the underlying probability, resp. expectation operator, will be denoted by $\mathbb{P}_\theta$, resp. $\mathbb{E}_\theta$.

For a given McQN, we define the $\Theta$-\emph{stability region} via the associated Markov process $\mathcal{X}$, as follows:
$$\Theta_{\rm s}:=\{\theta\in\Theta:\:\mathcal{X}\:\text{is stable under}\:\mathbb{P}_\theta\};$$
here, by \emph{stability} we mean positive (Harris) recurrence. Stability of McQNs has been thoroughly investigated in
\cite{Dai:95,BGT:96,H:97,GH:05,Bramson:08}.

\begin{remk} 
{\em Note that the concept ``stability region'' is slightly different from the one introduced in \cite{Dai}, which refers to the
stability of the associated fluid model; the latter is, in general, a subset of the former \cite{Dai:95}.}
\end{remk}

In the same vein, define the $\Theta$-\emph{subcriticality region}
$$\Theta_{\rm c}:=\left\{\theta\in\Theta:\:\max_{i}\sum_{k\in\mathcal{K}_i}\frac{[(I-R')^{-1}\theta]_k}{\beta_k}<1\right\}.$$
If $\Theta=\mathbb{R}^d_+$, we shall use the terminology \emph{full stability (subcriticality) region} and we shall omit specifying $\Theta$ when not relevant, or no confusion occurs.

In many cases (e.g. Jackson and Kelly networks) stability is equivalent to subcriticality, hence $\Theta_{\rm s}=\Theta_{\rm c}$.
Nevertheless, this is not always the case, as illustrated by numerous (counter) examples in the literature (see also Example 
\ref{Dai:exmp} below) and, in general, stability only implies subcriticality, hence $\Theta_{\rm s}\subseteq\Theta_{\rm c}$;
see \cite{Bramson:08}.

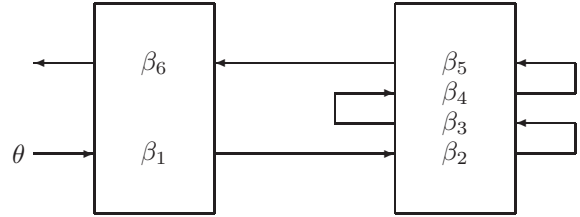
\begin{figure}
\unitlength=4mm \linethickness{1pt}
\begin{center}
\begin{picture}(20.0,08.0)
%SERVER 1
\drawline(03.0,0.0)(07.0,0.0)
\drawline(07.0,0.0)(07.0,7.0)
\drawline(03.0,7.0)(07.0,7.0)
\drawline(03.0,7.0)(03.0,0.0)
%SERVER 2
\drawline(13.0,0.0)(17.0,0.0)
\drawline(17.0,0.0)(17.0,7.0)
\drawline(13.0,7.0)(17.0,7.0)
\drawline(13.0,7.0)(13.0,0.0)
%ROUTING
\linethickness{0.4pt}
\put(01.0,2.0){\vector(1,0){2.00}}
\put(07.0,2.0){\vector(1,0){6.00}}
\drawline(17.0,2.0)(19.0,2.0)
\drawline(19.0,2.0)(19.0,3.0)
\put(19.0,3.0){\vector(-1,0){2.00}}
\drawline(13.0,3.0)(11.0,3.0)
\drawline(11.0,3.0)(11.0,4.0)
\put(11.0,4.0){\vector(1,0){2.00}}
\drawline(17.0,4.0)(19.0,4.0)
\drawline(19.0,4.0)(19.0,5.0)
\put(19.0,5.0){\vector(-1,0){2.00}}
\put(13.0,5.0){\vector(-1,0){6.00}}
\put(03.0,5.0){\vector(-1,0){2.00}}

%\dottedline(06.0,4.0)(10.0,4.0)
%\dottedline(06.0,7.0)(10.0,7.0)
%\dottedline(16.0,4.0)(20.0,4.0)
%\dottedline(16.0,5.0)(20.0,5.0)
%\dottedline(16.0,6.0)(20.0,6.0)
%\dottedline(16.0,7.0)(20.0,7.0)

%MARKS
%\put(08.0,1.0){\makebox(0,0)[cc]{$\mathfrak{S}_1$}}
%\put(18.0,1.0){\makebox(0,0)[cc]{$\mathfrak{S}_2$}}
\put(0.5,2.0){\makebox(0,0)[cc]{$\theta$}}
\put(05.0,2.0){\makebox(0,0)[cc]{$\beta_1$}}
\put(05.0,5.0){\makebox(0,0)[cc]{$\beta_6$}}
\put(15.0,2.0){\makebox(0,0)[cc]{$\beta_2$}}
\put(15.0,3.0){\makebox(0,0)[cc]{$\beta_3$}}
\put(15.0,4.0){\makebox(0,0)[cc]{$\beta_4$}}
\put(15.0,5.0){\makebox(0,0)[cc]{$\beta_5$}}
\end{picture}\end{center}
\caption{A first-come-first-serve reentrant line (Bramson -- Dai).}\label{BD:fig}
\end{figure}

\begin{exmp}\label{Dai:exmp}
{\em Consider a re-entrant line with two servers and six classes, with the routing indicated in Figure \ref{BD:fig}. Both stations employ the usual first-come-first-serve discipline.
We let $\theta\in\Theta=\{(r,0,\ldots,0):r\geq 0\}$, i.e.\ $r$ denotes the (Poisson) arrival rate, and denote by $\mu_1,\ldots,\mu_6$ the expected service times of the respective classes. Then, we have $$\rho_1=r(\mu_1+\mu_6),\:\rho_2=r(\mu_2+\mu_3+\mu_4+\mu_5).$$ However, if $\mu_1=\mu_3=\mu_4=\mu_5=0.001$,
$\mu_2=0.897$ and $\mu_6=0.899$, then $(1,0,\ldots,0)\in\Theta_{\rm c}\setminus\Theta_{\rm s}$, cf.~\cite{Dai:95}.}
\end{exmp}

We conclude that, except from the situations when $\Theta_{\rm s}=\Theta_{\rm c}$, no analytical representations are available,
in general, for stability regions. Therefore, numerical methods are sought instead. It is also worth noting that the full
subcriticality region is an open, bounded, star-shaped domain in $\mathbb{R}^d$, around the origin (vantage point); it is not
clear, however, whether the full stability region enjoys similar properties.
%that is, if $\theta\in\Theta_{\rm c}$ then the whole segment $[\mathbf{0},\theta]$ lies in $\Theta_{\rm c}$.

\subsection{Stability and Subcriticality Thresholds}

A vector $\vec{v}:=(v_1,\ldots,v_d)\in\mathbb{R}^d_+$ satisfying $\|\vec{v}\|=1$ will be called a (positive) \emph{direction} in $\mathbb{R}^d$; for a given direction $\vec{v}$,
we define the $\vec{v}$-\emph{ray} $$\langle\vec{v}\,\rangle:=\{r\cdot\vec{v}:r\geq 0\};$$ the $\vec{v}$-ray is a one-dimensional manifold isomorphic to $[0,\infty)$, hence one can endow it with the usual ordering and topology on the real non-negative half-line.

In the sequel, we shall restrict our analysis to the case $\Theta=\langle\vec{v}\,\rangle$; there are at least two reasons for that:
\begin{itemize}
\item The family of all rays $\langle\vec{v}\,\rangle$ sweeps the whole non-negative quadrant, hence any given set $\mathfrak{D}\subseteq\mathbb{R}^d_+$ is characterized by
the family of traces it leaves on the positive rays.
%; in particular, if $\mathfrak{D}$ is a bounded star-shaped
%domain, having the origin as vantage point, these traces are simply segments connecting the origin to boundary points and one can reconstruct the domain $\mathfrak{D}$ by interpolating such points.
\item Many of the practical applications of McQNs concern reentrant lines, where $\Theta=\langle(1,0,\ldots,0)\rangle$.
\end{itemize}

For an arbitrary positive direction $\vec{v}$, we define the \emph{stability threshold} in direction $\vec{v}$ as $\theta_*(\vec{v}\,):=\sup\langle\vec{v}\,\rangle_{\rm s}$;
in the same vein, we define the \emph{critical threshold} in in direction $\vec{v}$ as $\bar{\theta}(\vec{v}\,):=\sup\langle\vec{v}\,\rangle_{\rm c}$.

When the direction $\vec{v}$ is not relevant (or clear from the context), we shall use the simplified notations $\theta_*$, resp. $\bar{\theta}$; we stress however that both
thresholds depend on $\vec{v}$. Note that, in the light of the properties put forward in Section \ref{regio:sec}, it holds that $\mathbf{0}<\theta_*\leq\bar{\theta}<\infty$;
the leftmost inequality follows from the fact that the full stability region includes the open set
$$\left\{\theta\in\mathbb{R}^d_+:\:\sum_{k=1}^d\frac{[(I-R')^{-1}\theta]_k}{\beta_k}<1\right\},$$
corresponding to the sufficient (global) stability condition $\rho_1+\ldots+\rho_\aleph<1$; see \cite{Bramson:08}.

Furthermore, we have $\langle\vec{v}\,\rangle_{\rm c}=[\mathbf{0},\bar{\theta})$,
but a similar representation does not necessarily hold for $\langle\vec{v}\,\rangle_{\rm s}$, unless the full stability region has similar geometric properties as the subcriticality
region, i.e., it is an open, star-shaped domain; it holds, however, that $\langle\vec{v}\,\rangle_{\rm s}\subseteq[\mathbf{0},\theta_*)$.

Finally, we note that for any direction $\vec{v}$ there exist finite (positive) constants $r_*$ and $\bar{r}$ (both depending on $\vec{v}$), such that
$\theta_*(\vec{v}\,)=r_*(\vec{v}\,)\cdot\vec{v}$, resp.\ $\bar{\theta}(\vec{v}\,)=\bar{r}(\vec{v}\,)\cdot\vec{v}$; in addition, $\bar{r}$ can always be analytically calculated,
as follows: $\bar{r}(\vec{v}\,)=\min_i\bar{r}_i(\vec{v}\,)$, where
\begin{equation}\label{critic:eq}
\bar{r}_i(\vec{v}\,):=\left[\sum_{k\in\mathcal{K}_i}\frac{\delta_k}{\beta_k}\right]^{-1},
\end{equation}
denotes the critical threshold for station $i$; in the last display, we used the short-hand notation
$$\delta:=(I-R')^{-1}\vec{v}=(I+R+R^2+\ldots)'\vec{v};$$
in particular, $\theta=r\cdot\vec{v}$ entails $\lambda=r\delta$ on $\langle\vec{v}\,\rangle$.

\newpage

\section{A Numerical Method for Determining Stability Regions}\label{main:sec}

Throughout this section, $\mathcal{X}$ will denote the (Markov) queue-configuration process associated with an McQN with $\aleph$
stations, $d$ classes, arrival-rate vector $\theta$, service-rate vector $\beta$ and routing matrix $R$, while $\vec{v}$ will denote
a fixed positive direction in $\mathbb{R}^d$; in particular, $\theta=r\cdot\vec{v}$.

The aim is to design a numerical method for evaluating the stability threshold $\theta_*$ along the positive direction $\vec{v}$. Our
analysis will reveal that, under some (rather weak) monotonicity conditions, the $\langle\vec{v}\,\rangle$-stability region satisfies
$\langle\vec{v}\,\rangle_{\rm s}=[\mathbf{0},\theta_*)=[0,r_*)\cdot\vec{v}$ and that the stability threshold $\theta_*$ (in fact, $r_*$)
can be evaluated via Robbins-Monro schemes; eventually, we extend this method to more general, star-convex parameter sets.

\subsection{Stability Thresholds for Jackson Networks}\label{Jackson:sec}

Assume that $\mathcal{X}$ corresponds to a Jackson network with $d$ stations/classes. In this case, $\Theta_{\rm s}=\Theta_{\rm c}$,
for any $\Theta$, hence $\theta_*=\bar{\theta}$, resp. $r_*=\bar{r}=\min_k(\beta_k/\delta_k)$, cf.~(\ref{critic:eq}).
Furthermore, consider $\phi:\mathbb{X}=\mathbb{N}^d\longrightarrow(0,1]$ defined as $\phi(\mathbf{x}):=\exp(-\alpha\|\mathbf{x}\|)$,
for some (fixed) $\alpha>0$; then for any $r<\bar{r}$ (stability) it holds (cf.\ \cite{Jackson}) that (recall that $\theta=r\cdot\vec{v}$)
\begin{equation}\label{equilibrium:eq}
\tilde{\varphi}(r):=\lim_{t\rightarrow\infty}\mathbb{E}_\theta[\phi(X_t)]=\prod_{k=1}^d\frac{\bar{r}_k-r}{\bar{r}_k-re^{-\alpha}}.
\end{equation}
The function $\tilde{\varphi}$ in the above display is continuous and strictly decreasing on $[0,\bar{r})$, with
$\tilde{\varphi}(0)=1$, $\tilde{\varphi}(\bar{r})=0$.\\ Therefore, denoting by
$r_\varepsilon\in(0,\bar{r})$ the (unique) root of the equation $\tilde{\varphi}(r)=\varepsilon\in(0,1)$, we note that
$r_\varepsilon$ is increasing in $\varepsilon$ and it can be verified that it approaches $\bar{r}$ as $\varepsilon$ decreases to $0$;
the same holds true if we replace $\phi$ by any bounded function vanishing at infinity.

One concludes that, for Jackson networks, stability thresholds can be approximated by roots of equations of the type $\tilde{\varphi}(r)=\varepsilon$
(for $\varepsilon$ close to $0$), where $\tilde{\varphi}$ is a stationary performance measure of the network under consideration;
more specifically,
$\tilde{\varphi}$ appears as the expectation under the equilibrium distribution of some bounded function vanishing at infinity.

\subsection{Stability Thresholds in the Multi-class Setup}

In this section we extend the approximation scheme described in Section \ref{Jackson:sec} beyond the Jackson network setup, in order to approximate
stability thresholds in cases where they are not available in closed form. In doing so, the following questions/challenges arise:
\begin{enumerate}
\item [(I)] Is the $\langle\vec{v}\,\rangle$-stability region still a half-open interval of the form $[\mathbf{0},\theta_*)=[0,r_*)\cdot\vec{v}$, so that
$r_*$ determines the stability region?
\item [(II)] Provided   the answer in (I) is affirmative, does there exist a (stationary) performance measure $\tilde{\varphi}:[0,\infty)\longrightarrow[0,1]$
such that the root $r_\varepsilon$ of the equation $\tilde{\varphi}(r)=\varepsilon$ approaches the threshold $r_*$, for $\varepsilon$ close to $0$?
\item [(III)] Provided  the answers in (I)--(II) are affirmative, how to evaluate $r_\varepsilon$, since analytical expressions for stationary performance measures,
such as the one in (\ref{equilibrium:eq}), are not available in general?
\end{enumerate}
In what follows, we shall provide a set of conditions guaranteeing positive answers to questions (I) and (II) above and discuss possible approaches to (III).

To start with, note that (I) assumes a certain type of monotonic behavior. More specifically, it requires that the stability region
is a monotone set, in the sense that stability for a certain parameter entails stability for all ``smaller'' parameters. In addition,
if the answer is affirmative for any direction $\vec{v}$, then the full stability region defines a star-shaped domain around the origin.

%Below we assume the componentwise ordering on $\Theta$; when $\Theta=\langle\vec{v}\rangle$ this reduces to the ordering induced by $r\mapsto r\cdot\vec{v}$.

Assume now that there exists some $\phi:\mathbb{X}\longrightarrow(0,1]$, vanishing at infinity, satisfying the following condition:
\begin{enumerate}
\item [(\textbf{M}1)] the mapping
$$(t,\theta)\longmapsto\varphi_t(\theta):=\mathbb{E}_\theta[\phi(X_t)|X_0=\emptyset],$$ is (jointly) non-increasing on $[0,\infty)\times\Theta$;
\end{enumerate}
that is, we assume the existence of some functional $\phi$ of the process $\mathcal{X}$ (started in the empty configuration) which is monotone (in expectation)
w.r.t.~both time and arrival-rates (componentwise ordering).

Provided that (\textbf{M}1) above holds true, the limit
\begin{equation}\label{limit:eq}
\varphi(\theta):=\lim_{t\rightarrow\infty}\varphi_t(\theta)=\inf_{t\geq 0}\varphi_t(\theta)\in[0,1],
\end{equation}
exists and defines a non-decreasing function on $\Theta$. For any such $\varphi$ it holds that
\begin{equation}\label{regio:eq}
\Theta_{\rm s}=\{\theta\in\Theta:\varphi(\theta)>0\};
%\langle\vec{v}\rangle_{\rm s}=\{r\geq 0:\tilde{\varphi}(r)>0\}=[0,r_*),
\end{equation}
in particular, given the positive direction $\vec{v}$, we define
$\tilde{\varphi}_t,\tilde{\varphi}:[0,\infty)\longrightarrow[0,1]$ as the push-forwards of $\varphi_t$, resp. $\varphi$, on the ray
$\vec{v}$; that is, $\tilde{\varphi}(r)=\varphi(r\cdot\vec{v})$. Then,
$$\langle\vec{v}\rangle_{\rm s}=\tilde{\varphi}^{-1}((0,1])\cdot\vec{v}=[0,r_*)\cdot\vec{v},$$
which solves question (I). A complete proof of the above facts is provided in \cite{LM:2016}.
% and relies on the fact that a Markov process on discrete state-space is either stable or leaves eventually any compact set, almost surely. %\cite{Meyn-Tweedie:II}.

\newpage

Furthermore, to guarantee (II), it suffices that
\begin{enumerate}
\item [(\textbf{M}2)] the mapping $\varphi:\Theta\longrightarrow[0,1]$ defined by (\ref{limit:eq}) is continuous and strictly
decreasing on $\Theta_{\rm s}$.
\end{enumerate}
Indeed, assuming that (\textbf{M}2) holds true, the function $\tilde{\varphi}:[0,r_*)\longrightarrow(0,1]$ is homeomorphic, hence
the root $r_\varepsilon=\tilde{\varphi}^{-1}(\varepsilon)$ is correctly defined and approximates the threshold $r_*$, for
$\varepsilon\rightarrow 0$.

Finally, for estimating the root $r_\varepsilon$ one can employ a stochastic approximation scheme of Robbins-Monro (RM) type \cite{RM},
which requires that the values of $\varphi$ (for various parameters) are evaluated by simulation. More specifically, an RM approximation
scheme is an iterative method which constructs a sequence of parameter updates such that at every update an unbiased estimate of
$\varphi$ is used to generate a new parameter. The main difficulty when applying an RM scheme in this setting arises from the fact
that one needs to sample from $\varphi$, which appears as a stationary (limiting) measure of the process $\mathcal{X}$. There are
two possible approaches:
\begin{enumerate}
\item[(1)] direct simulation via regenerative ratios, which in turn requires simulating the queue-configuration process along a
regenerative cycle; see e.g.\ \cite{Asm-Glynn}.
\item[(2)] simulating instead $\varphi_t$, for an increasing sequence of time-horizons $t\rightarrow\infty$ and invoking an approximation
argument; see e.g.\ \cite {Burk:56}.
\end{enumerate}
Method (1) seems more forthright. Note however that recurrence times are random and may become arbitrarily large as the input
parameter  approaches the boundary of the stability region. Since one expects that the approximation scheme will stabilize
somewhere in the neighborhood  of the stability threshold, i.e. at the boundary of the stability region, such a method
seems rather unpredictable in terms of computational effort. Method (2) avoids this inconvenience by setting fixed simulation
horizons, hence allows for a better control over the computational complexity.

We conclude this section with several considerations on the two conditions formulated above:% which are in order here.
\begin{itemize}
\item Conditions (\textbf{M}1) and (\textbf{M}2) are deliberately stated for general $\Theta$ (rather than $\langle\vec{v}\,\rangle$),
since in many situations, the conditions hold for $\Theta=\mathbb{R}^d_+$, which entails their validity (for the same $\phi$) for any ray.
\item Conditions (\textbf{M}1) and (\textbf{M}2) are quite common for Jackson networks; (\textbf{M}1) follows by standard stochastic
monotonicity theory for Markov chains, whereas (\textbf{M}2) follows directly by (\ref{equilibrium:eq}).
\item Thm.\ 1 in \cite{LM:2016} establishes the validity of (\textbf{M}1), provided that the queue-configuration process $\mathcal{X}$
fulfils a certain stochastic monotonicity condition.
\item Prop.\ 1 and 2 in \cite{LM:2016} show that, for a wide class of McQNs (including the examples treated in this paper),
conditions (\textbf{M}1) and (\textbf{M}2) hold for certain $\phi$'s (hence, $\varphi$'s) and for $\Theta=\mathbb{R}^d_+$.
\end{itemize}

\subsection{Numerical Evaluation of Stability Thresholds}\label{method:sec}

In this section we assume that conditions (\textbf{M}1) and (\textbf{M}2) hold for $\Theta=\langle\vec{v}\rangle$, for a certain
$\phi:\mathbb{X}\longrightarrow(0,1]$, vanishing at infinity and we design a numerical method for approximating the stability
threshold $\theta_*=r_*\cdot\vec{v}$.%; this allows one to recover the stability region $[\mathbf{0},\theta_*)$.

Fix some arbitrary increasing sequence $\{t_n\}_{n\geq 0}$ of non-negative numbers satisfying $t_0=0$, $t_n\rightarrow\infty$ and let
$\mathfrak{D}_n(r)$ denote the distribution of $\phi(X_{t_n})$ under $\mathbb{P}_\theta$, for $\theta=r\cdot\vec{v}$ and $n\geq 0$;
we further set $\mathfrak{D}_n(r)=\mathfrak{D}_n(0)$, for $r<0$.
Furthermore, fix some sequence $\{a_n\}_{n\geq 1}$ of decreasing positive numbers satisfying $a_n\rightarrow 0$ and $\varepsilon>0$
and define the sequence of iterates %(for $n\geq 0$)
\begin{equation}\label{st-approx:eq}
\forall n\geq 1:\: x_n = x_{n-1} + a_n\left(z_n-\varepsilon\right),
\end{equation}
%\begin{equation}\label{st-approx:eq}
%x_{n+1}=\min\{\bar{r},\max\left\{x_n+a_{n+1}\left(z_n-\varepsilon\right),0\right\}\},
%\end{equation}
where $x_0\in(0,\bar{r})$ is arbitrarily chosen and for each $n\geq 0$ the r.v.\ $z_n$ follows the conditional distribution
$\mathfrak{L}[z_n|x_{n-1}]=\mathfrak{D}_n(x_{n-1})$, given $x_{n-1}$. 
%Note that, $x_n\in[0,\bar{r}]$, for $n\geq 0$, by construction.

Our next result establishes the convergence of the iterates $x_n$ in (\ref{st-approx:eq}) towards the root $r_\varepsilon$ of
the equation $\tilde{\varphi}(r)=\varphi(r\cdot\vec{v})=\varepsilon$, for $n\rightarrow\infty$ and, under slightly more restrictive
conditions, provides the magnitude of the approximation error; see the Appendix for a proof.

\begin{them}\label{conv:them}
For any $\varepsilon>0$, $x_0\in(0,\bar{r})$ and positive sequences $\{t_n\}_{n\geq 0}$ and $\{a_n\}_{n\geq 1}$, satisfying
$$\lim_{n\rightarrow\infty}t_n=\infty,\:\sum_{n\geq 1}a_n=\infty,\:\sum_{n\geq 1}a_n^2<\infty,$$
the iterates $\{x_n\}_{n\geq 0}$ in (\ref{st-approx:eq}) satisfy $x_n\longrightarrow r_\varepsilon$, a.s.

Furthermore, assume that the family of derivatives $\tilde{\varphi}_t'$ converges uniformly on $(0,r_*)$, for $t\rightarrow\infty$,
and that $\inf|\tilde{\varphi}'(r)|>0$. If $a_n=a\cdot n^{-\omega}$, for $\omega\in(1/2,1]$ and $a>(\omega-1/2)/(\inf|\tilde{\varphi}_t'(r)|)$ and 
\begin{equation}\label{conv:eq}
\sup_{r\in[0,r_*]}|\tilde{\varphi}_{t_n}(r)-\tilde{\varphi}(r)|=o(n^{-\kappa}),
\end{equation} 
for some $\kappa>\omega-1/2$, then (in probability) $$(x_n-r_\varepsilon)=O(n^{-(\omega-1/2)}).$$
\end{them}

%We conclude this section by discussing the approximation error of the proposed method. Theorem \ref{conv:them}
In practice, we fix some large $n\geq 1$ and use the estimate $x_n$ to approximate $r_*$. The approximation error consists
of a random and a deterministic component:
\begin{equation}\label{error:eq}
\Delta_n^{\varepsilon}:=|x_n-r_*|\leq|x_n-r_\varepsilon|+(r_*-r_\varepsilon).
\end{equation}
While the behavior of the random component $|x_n-r_\varepsilon|$ is established by Theorem \ref{conv:them}, for the deterministic
part in (\ref{error:eq}) we note that (for small $\varepsilon$)
$$\varepsilon=\tilde{\varphi}(r_\varepsilon)-\tilde{\varphi}(r_*)\approx -\tilde{\varphi}'(r_\varepsilon)(r_*-r_\varepsilon);$$
in particular, if $\tilde{\varphi}'$ is bounded away from $0$ (close to $r_*$) then one obtains $r_\varepsilon\rightarrow r_*$ at
a linear rate. Nevertheless, if $\lim_{r\rightarrow r_*}\tilde{\varphi}'(r)=0$ then convergence is slower, as we shall note in our
numerical experiments in Section \ref{num:sec}.

We conclude that the approximation error of the method depends essentially on the behavior of the derivative $\tilde{\varphi}'$ close
to $r_*$; more specifically, denoting $c:=\lim_{r\rightarrow r_*}|\tilde{\varphi}'(r)|$, we note that the larger $c$, the better the accuracy.

\subsection{Approximating Star-shaped Stability Regions}

In this section we assume that conditions (\textbf{M}1) and (\textbf{M}2) hold true for a certain $\phi:\mathbb{X}\longrightarrow(0,1]$
and for some star-shaped (around the origin) parameter set $\Theta\subseteq\mathbb{R}^d_+$. Then the $\Theta$-stability region
$\Theta_{\rm s}$ defines itself a star-shaped domain around the origin; a similar fact holds for the stability region for the fluid model
\cite{Chen}. 

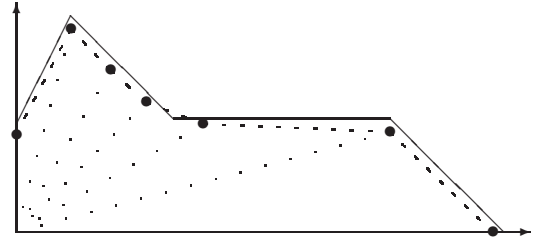
\begin{figure}
\unitlength=3.6mm
\begin{center}
\begin{picture}(19.0,8.5)
\linethickness{0.6pt}
%AXES
\put(0.0,0.0){\vector(1,0){19}}
\put(0.0,0.0){\vector(0,1){8.5}}
%RAYS
\dottedline(0,0)(2.02,7.5)%(2.7,10.0)
\dottedline(0,0)(3.48,6.0)%(5.8,10.0)
\dottedline(0,0)(4.8,4.8)%(10,10.0)
\dottedline(0,0)(6.92,4.0)%(17.3,10.0)
\dottedline(0,0)(13.8,3.7)%(20.0,5.36)
%MARKS
\put(0.0,3.6){\makebox(0,0)[cc]{$\bullet$}}
\put(2.02,7.5){\makebox(0,0)[cc]{$\bullet$}}
\put(3.48,6.0){\makebox(0,0)[cc]{$\bullet$}}
\put(4.8,4.8){\makebox(0,0)[cc]{$\bullet$}}
\put(6.9,4.0){\makebox(0,0)[cc]{$\bullet$}}
\put(13.8,3.7){\makebox(0,0)[cc]{$\bullet$}}
\put(17.6,0.0){\makebox(0,0)[cc]{$\bullet$}}
%INTERPOLATION
\dashline{0.2}(0.0,3.6)(2.02,7.5)
\dashline{0.2}(2.02,7.5)(3.48,6.0)
\dashline{0.2}(3.48,6.0)(4.8,4.8)
\dashline{0.2}(4.8,4.8)(6.92,4.0)
\dashline{0.2}(6.92,4.0)(13.8,3.7)
\dashline{0.2}(13.8,3.7)(17.6,0.0)
%BOUNDARY
\linethickness{0.2pt}
\drawline(0.0,4.0)(2.0,8.0)
\drawline(2.0,8.0)(5.8,4.2)
\drawline(5.8,4.2)(13.8,4.2)
\drawline(13.8,4.2)(18.0,0.0)
\end{picture}\end{center}
\caption{Numerical approximation of a star-shaped domain. Solid line = true boundary; dotted lines = rays; bullets = approximations
of the boundary points (thresholds); dashed line = approximated boundary.}\label{star:fig}
\end{figure}

Assuming w.l.o.g.\ that $\Theta=\mathbb{R}_+^{d'}$, for some $d'\leq d$, such a domain can be approximated as follows: one can construct 
a grid of points on the positive orthant of the unit sphere (each point corresponding to a given direction) and determine the stability threshold along each direction, cf.~Section \ref{method:sec}. Finally, one connects thresholds corresponding to neighboring points (directions), obtaining in this way a polytope which approximates the $\Theta$-stability region (for large number of points); such a procedure, for $d'=2$, is graphically illustrated in Figure \ref{star:fig}.

\begin{remk} {\em Note that the boundary point of the $\Theta$-stability region in some given direction is obtained as the minimum
between the boundary point of $\Theta$ and the stability threshold in that direction; for a more efficient numerical procedure, one
can replace $\bar{r}$ in (\ref{st-approx:eq}) by the corresponding boundary point of $\Theta$, thus avoiding to simulate (too) congested networks.}
\end{remk}

\vspace{-4mm}

\begin{figure}[h!]
\unitlength=4mm \linethickness{1pt}
\begin{center}
\begin{picture}(20.0,8.0)
%SERVER 1
\drawline(7.0,0.0)(3.0,0.0)\drawline(7.0,0.0)(7.0,7.0)
\drawline(3.0,7.0)(7.0,7.0)\drawline(3.0,7.0)(3.0,0.0)
%SERVER 2
\drawline(13.0,0.0)(17.0,0.0)\drawline(17.0,0.0)(17.0,7.0)
\drawline(13.0,7.0)(17.0,7.0)\drawline(13.0,7.0)(13.0,0.0)
%ROUTING
\linethickness{0.4pt}
\put(1.0,2.0){\vector(1,0){2.00}}
\put(7.0,2.0){\vector(1,0){6.00}}
\drawline(17.0,2.0)(19.0,2.0)\drawline(19.0,2.0)(19.0,5.0)
\put(19.0,5.0){\vector(-1,0){2.00}}
\put(13.0,5.0){\vector(-1,0){6.00}}
\put(3.0,5.0){\vector(-1,0){2.00}}

%MARKS
\put(0.5,2.0){\makebox(0,0)[cc]{$\theta$}}
\put(3.6,2.0){\makebox(0,0)[cc]{$\nabla$}}
\put(5.0,2.0){\makebox(0,0)[cc]{$\beta_1$}}
\put(6.4,5.0){\makebox(0,0)[cc]{$\Diamond$}}
\put(5.0,5.0){\makebox(0,0)[cc]{$\beta_4$}}
\put(13.6,2.0){\makebox(0,0)[cc]{$\Diamond$}}
\put(15.0,2.0){\makebox(0,0)[cc]{$\beta_2$}}
\put(16.4,5.0){\makebox(0,0)[cc]{$\nabla$}}
\put(15.0,5.0){\makebox(0,0)[cc]{$\beta_3$}}
\end{picture}
\end{center}
\caption{\label{LK:fig} The Lu-Kumar network; $\Diamond$ gives priority over $\nabla$.}
\end{figure}
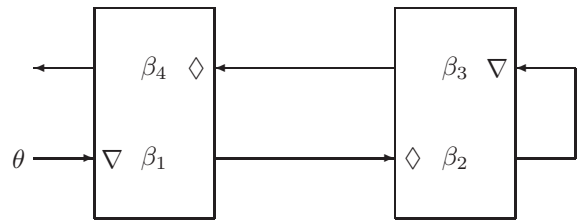

\newpage

\section{Numerical Results}\label{num:sec}

%In this section we illustrate the use of the method developed in Section \ref{method:sec}. We include here results
%which correspond to a set of meaningful examples; that is, we consider both trivial examples (for which the stability
%thresholds are known) and cases in which stability thresholds are not known. The algorithm parameters are chosen such
%that satisfactory approximations are obtained in the trivial cases; however, drawing conclusions in the non-trivial cases
%is conditioned on the assumption that condition (\ref{conv:eq}) (Them.~\ref{conv:them}) holds true.

In this section we illustrate the use of the method developed in Section \ref{method:sec}. We include here experimental
results corresponding to examples for which the stability thresholds are known, and results for which these are not known.

For a given $\varepsilon>0$ we average out $N=10000$ RM iterates (\ref{st-approx:eq}) in order to construct an estimator 
\begin{equation}\label{estimator:eq}
\hat{r}_\varepsilon:=\frac{1}{N}\sum_{n=1}^N x_n,
\end{equation}
for the solution $r_\varepsilon$ of $\varphi(r)=\varepsilon$, where $\varphi$ is defined by (\ref{limit:eq}), for $\phi(\xi)=\exp(-\alpha\|\xi\|)$, with $\|\xi\|$ denoting the total number of jobs in the network configuration
$\xi$; the average in (\ref{estimator:eq}) has the advantage that it is less sensitive to initial jumps/outliers
\cite{PJ:92}. For illustrative purposes, we analyze the effect of varying the value of $\varepsilon$; more specifically,
we let $\varepsilon=10^{-c}$, with $c=7,8,9,10$.

For the numerical experiments below, $a_n=a/n^\omega$ and $t_n=t_0+bn$, with  $\omega=1$, $a=1/\varepsilon$, $t_0=2\cdot 10^6$,
$b=200$; also, let $x_0=0$ and $\alpha=1$. These parameters are set such that they provide (approximately) correct values when 
the stability region is known; note that otherwise firm conclusions can only be drawn under the proviso that condition (\ref{conv:eq}) 
(Thm.\ \ref{conv:them}) holds true.

\begin{table}
\begin{center}
\begin{tabular}{||c||c|c|c||c||}
\hline
$\log\varepsilon$    & $-2$ & $-4$ & $-6$ & $\bar{r}=r_*$ \\
\hline\hline
$v=0.000$ & $1.9730$ & $1.9984$ & $1.9999$ & $2.0000$ \\
\hline
$v=0.268$ & $1.9809$ & $1.9984$ & $1.9998$ & $2.0000$ \\
\hline
$v=0.577$ & $1.9004$ & $1.9900$ & $1.9992$ & $2.0000$ \\
\hline
$v=1.000$ & $1.3132$ & $1.3319$ & $1.3332$ & $1.3333$ \\
\hline
$v=1.732$ & $0.8178$ & $0.8271$ & $0.8280$ & $0.8282$ \\
\hline
$v=3.732$ & $0.4013$ & $0.4061$ & $0.4068$ & $0.4069$ \\
\hline
\end{tabular}
\end{center}

\vspace{-4mm}

\caption{\label{Jackson:tab}\footnotesize{Stability threshold estimates along various slopes ($v$), for various accuracy levels
($\varepsilon$) for the Jackson network}.}

\end{table}

%All the numerical results presented in this section are based on simulations performed in parallel on the DAS-4 architecture: 
%\url{http://www.cs.vu.nl/das4}.

\subsection{Jackson Networks}

Consider an open Jackson network consisting of two servers/classes $k=1,2$, having input rates $\theta_1$, resp.~$\theta_2$, and
service rates $\beta_1$, resp.~$\beta_2$.
We further assume that any job finishing service at server $1$ moves to server $2$ with probability $\wp\in[0,1]$, or leaves the network; that is, $R_{12}=\wp$ and $R_{11}=R_{21}=R_{22}=0$.

Pick now some $\vec{v}=(1,v)$, for some arbitrary $v\geq 0$ and recall that $r_*=\bar{r}=\bar{r}_1\wedge\bar{r}_2$, where,
cf.~(\ref{critic:eq}),
\begin{equation}\label{threshold:eq}
\bar{r}_1=\beta_1\cdot\|\vec{v}\|,\:\bar{r}_2=(\wp+v)^{-1}\beta_2\cdot\|\vec{v}\|.
\end{equation}

Let $\beta_1=2$, $\beta_2=1.6$, $\wp=0.2$. A summary of the corresponding results compared to the true values, calculated
using (\ref{threshold:eq}), is provided in Table \ref{Jackson:tab}.

\subsection{Multi-class Reentrant Lines}

For reentrant lines, $\Theta=\langle\vec{v}\rangle$, with $\vec{v}=(1,0,\ldots,0)$, so that under the monotonicity condition
(\textbf{M}1) the $\Theta$-stability region is determined (only) by the stability threshold $r_*$. For the networks considered
below, it has been demonstrated in \cite{LM:2016} that monotonicity conditions (\textbf{M}1) and (\textbf{M}2) hold; however,
for illustrative purposes, here we test condition (\textbf{M}1) numerically; see Table \ref{allTheta}.

Our first example is the network in Example \ref{Dai:exmp}, for which $\bar{r}=\bar{r}_1=\bar{r}_2=1/0.9\simeq 1.111$,
cf.~(\ref{critic:eq}). Table \ref{allR} (A) displays estimates $\hat{r}_\varepsilon$ for the above specified $\varepsilon$'s.

Secondly, consider the Lu-Kumar network \cite{LK:91}, in which both stations employ a (preemptive) priority policy, as illustrated
in Figure \ref{LK:fig}. Stability holds iff the network is subcritical and $\theta(\beta_2+\beta_4)<\beta_2\beta_4$ \cite{Dai:96}.
For our numerical experiments, we let $\beta_1=1.2$, $\beta_3=2$ and $\beta_2=\beta_4=1$, hence $\bar{r}=0.545$ and $r_*=0.5$;
this is illustrated in Table \ref{allR} (B).

Finally, consider the FCFS version of the Lu-Kumar network, with the same service rates; in this case, determining the stability
region is an open problem, cf.~\cite{Dai:96}. Our numerical results, provided in Table \ref{allR} (C), suggest that stability and subcriticality are equivalent.

The estimates $\hat{r}_\varepsilon$ in Table \ref{allR} provide approximations for $r_\varepsilon(t)$, the root of the equation $\varphi_t(r)=\varepsilon$, where $t=t_N=t_0+bN$, which in turn approximate 
$r_\varepsilon$ (for large $t$). Furthermore, it holds that $$\lim_{\varepsilon\rightarrow 0} r_\varepsilon(t)=\bar{r},\quad 
\lim_{\varepsilon\rightarrow 0}\lim_{t\rightarrow\infty}r_\varepsilon(t)=r_*.$$
In particular, the limits above are not interchangeable when $r_*\neq\bar{r}$ and the iterates $\hat{r}_\varepsilon$ do not 
converge to $r_*$ in these cases, as suggested by Table \ref{allR} (A) and (B).

\section{Concluding Remarks}\label{remk:sec}

In this paper we have developed a simulation-based numerical method for determining the stability region (w.r.t.~arrival rates) associated with Markovian McQNs. Our method identifies thresholds at which the queue sizes `explode'. In particular, stability
regions for networks for which no analytical stability conditions are known can be approximated numerically. The method does not 
extend in a straightforward way to the non-Markovian McQNs (non-exponential distributions), as the required (stochastic) monotonicity
properties for such networks have not been not established yet.

The complexity of a given network is reflected by the number of stations, classes and positive entries in the routing matrix. The computation time for generating one iterate $x_n$ increases linearly w.r.t.\ the time-horizon $t_n$. The trade-off between method complexity and accuracy is governed by the growth rate of the sequence $\{t_n\}_n$, hence gaining insight into the impact of
the choice of $\{t_n\}_n$ deserves future research efforts.

\appendix

\noindent\textbf{Proof} of Thm.\ \ref{conv:them}: The proof is based on Thms.\ 1 and 2 in \cite{Burk:56}. Namely, for every $n\geq 0$ and $r\geq 0$, let us define the mean, resp. the variance:
$$\omega_n(r):=\mathbb{E}\left[\left(\varepsilon-z_n\right)|r\right],\:\sigma_n(r):=\mathrm{Var}\left[\left(\varepsilon-z_n\right)|r\right].$$
By the monotonicity assumption (\textbf{M}1), $\omega_n$ is continuous and increasing w.r.t.~$r\in(0,\bar{r})$ and $n$.
In addition, $\sigma_n(r)\leq 1$, for any $n$ and $r\geq 0$.

For the convergence part we apply Thm.\ 1 in \cite{Burk:56}; to this end, we verify the following set of conditions:
\begin{enumerate}
\item [(i)] $\omega_n,\sigma_n:[0,\infty)\longrightarrow\mathbb{R}$ are measurable, s.t.
$$\sup_{(n,r)}\frac{|\omega_n(r)|}{1+r}<\infty,\:\sup_{(n,r)}\sigma_n(r)<\infty;$$
\item [(ii)] for any $\epsilon>0$ there exists $n_\epsilon\geq 1$ s.t.~$|r-r_\varepsilon|>\epsilon$ entails $(r-r_\varepsilon)\omega_n(r)>0$, for $n\geq n_\epsilon$;
\item [(iii)] $\sum_n a_n^2<\infty$ and for $0<\epsilon_1<\epsilon_2$ it holds that
$$\sum_{n\geq 1}a_n\left(\inf_{\epsilon_1<|r-r_\varepsilon|<\epsilon_2}|\omega_n(r)|\right)=\infty.$$
\end{enumerate}
Condition (i) is immediate since $\omega_n(r)\in(-1,1)$ and $\sigma_n(r)\in[0,1]$, for any $(n,r)$. Set
$\omega(r):=\varepsilon-\tilde{\varphi}(r)$ and note that $r_\varepsilon$ appears as the (unique) root of the equation
$\omega(r)=0$, with $\omega$ being (strictly) increasing in $r_\varepsilon$, cf. (\textbf{M}2). Let $\epsilon>0$; since $\omega(r_\varepsilon+\epsilon)>0$ and $\omega_n(r)\uparrow\omega(r)$,
for $n\rightarrow\infty$, it follows that there exists some $n_\epsilon\geq 1$ such that $n\geq n_\epsilon$ entails $\omega_n(r_\varepsilon+\epsilon)>0$, hence for any $r>r_\varepsilon+\epsilon$
it holds that $$(r-r_\varepsilon)\omega_n(r)\geq(r-r_\varepsilon)\omega_n(r_\varepsilon+\epsilon)>0.$$
On the other hand, $r<r_\varepsilon-\epsilon$ entails $$\omega_n(r)\leq\omega(r)<\omega(r_\varepsilon)=0,$$ for any $n$,
hence (ii) follows true, as well.

Finally, to verify (iii) we let $\epsilon_1<\epsilon_2$ and (as before) we choose $n_1\geq 1$ (depending only on $\epsilon_1$),
such that $\omega_n(r)>0$ for $r>r_\varepsilon+\epsilon_1$ and every $n\geq n_1$. Since $\omega_n(r)<0$ for $r\leq r_\varepsilon-\epsilon_1$ and $n\geq 1$, one obtains for $n\geq n_1$ and
$\epsilon_1<|r-r_\varepsilon|<\epsilon_2$
\begin{eqnarray*}
|\omega_n(r)|
%\min\left\{\inf_{\theta-\theta_\varepsilon\in(\delta_1,\delta_2)}U_n(\theta),\:\inf_{\theta_\varepsilon-\theta\in(\delta_1,\delta_2)}-U_n(\theta)\right\}
& = & \min\left\{\omega_n(r_\varepsilon+\epsilon_1),-\omega_n(r_\varepsilon-\epsilon_2)\right\} \\
& \geq & \min\{\varepsilon-\tilde{\varphi}_{t_n}(r_\varepsilon+\epsilon_1),\tilde{\varphi}(r_\varepsilon-\epsilon_2)-\varepsilon\};
\end{eqnarray*}
using $\lim_n\min\{u_n,v\}=\min\{\lim_n u_n,v\}$ yields
$$\inf_r|\omega_n(r)|\geq
\min\{\varepsilon-\tilde{\varphi}(r_\varepsilon+\epsilon_1),\tilde{\varphi}(r_\varepsilon-\epsilon_2)-\varepsilon\}>0,$$
where the infimum is taken w.r.t.~$\epsilon_1<|r-r_\varepsilon|<\epsilon_2$.
Hence, (iii) holds true, provided that $$\sum_{n\geq 1} a_n=\infty,\quad\sum_{n\geq 1} a_n^2<\infty;$$ this proves the first claim.

\newpage

For the second part, we invoke Thm.\ 2 in \cite{Burk:56}; to this end, we verify the following set of conditions:
\begin{enumerate}
\item [(i)] For any $n$, $\omega_n(r)$ is strictly increasing in $r$; in particular, there exists the root $r_{\varepsilon}^n$
of $\omega_n(r)=0$.
\item [(ii)] The function sequence $\{\tau_n\}_{n\geq 0}$, defined as
$$\tau_n(r):=\left\{
       \begin{array}{ll}
        (r-r_{n,\varepsilon})^{-1}\omega_n(r), & \hbox{$r\neq r_{n,\varepsilon}$;} \\
         -\varphi'(r_\varepsilon), & \hbox{$r=r_{n,\varepsilon}$,}
          \end{array}
           \right.$$
satisfies $\tau_n(r)\in[M_1,M_2]$, for all $n,r$, with $M_1>0$ and $\tau_n(x_n)\rightarrow-\tilde{\varphi}'(r_\varepsilon)$ for $x_n\rightarrow r_\varepsilon$.
\item [(iii)] There exists constants $0\leq M_3<M_4$ such that $M_3\leq\sigma_n(r)=\mathrm{Var}[z_n|r]\leq M_4$, for all $n,r$,
and $x_n\rightarrow r_\varepsilon$ entails $\sigma_n(x_n)\rightarrow\sigma>0$.
\item [(iv)] there exist $\kappa,\omega$, s.t.\ $(\omega-1/2)\in(0,\kappa)$ and
$$(r_{\varepsilon}^n-r_\varepsilon)=o(n^{-\kappa}),\:n^\omega a_n\rightarrow a>(\omega-1/2)/M_1.$$
\end{enumerate}

Condition (i) is immediate since $\omega_n(r)=\varepsilon-\tilde{\varphi}_{t_n}(r)$ and $\tilde{\varphi}_t$ decreases, with $\tilde{\varphi}(0)=1$,
vanishing at infinity.

To verify (ii), we note that since $\tilde{\varphi}_t'$ is continuous and non-vanishing on $[0,r_*]$, hence it is bounded away
from both infinity and $0$, for any $t>0$; moreover, since $\tilde{\varphi}_{t_n}'$ converges uniformly to $\tilde{\varphi}_t'$,
which is continuous, non-vanishing on $(0,r_*)$, it follows that $\tilde{\varphi}_{t_n}'$ is uniformly bounded away from both $0$
and infinity. Furthermore, if $x_n\rightarrow r_\varepsilon$, such that $x_n\neq r_{\varepsilon}^n$, for all $n$, we obtain
(mean value) $\tau_n(x_n)=-\tilde{\varphi}_{t_n}'(u_n)$, for some $u_n$ satisfying $|u_n-r_\varepsilon|<\epsilon$, for some
(small) $\epsilon>0$. The convergence $\tilde{\varphi}_{t_n}'(u_n)\rightarrow\tilde{\varphi}'(r_\varepsilon)$ follows from the
uniform convergence of the derivatives; the convergence is not affected if $x_n=r_{\varepsilon}^n$, for some $n$'s.

Furthermore, the variance converges uniformly, viz.
$$\sigma_n(r)\rightarrow\sigma(r):=\left\{
                                     \begin{array}{ll}
                                       \mathrm{Var}_{\pi_\theta}[\phi(X)], & \hbox{$r<r_*$;} \\
                                       0, & \hbox{$r\geq r_*$,}
                                     \end{array}
                                   \right.$$
where (recall) $\theta=r\cdot\vec{v}$ and $\pi_\theta$ denotes the equilibrium distribution under $\mathbb{P}_\theta$, for $\theta\in\langle\vec{v}\,\rangle_{\rm s}$.
%$$\sigma_n(r)=\mathbb{E}_{r\cdot\vec{v}}^\emptyset\left[\phi^2(X_{t_n})\right]-\mathbb{E}_{r\cdot\vec{v}}^\emptyset\left[\phi(X_{t_n})\right]^2\longrightarrow\sigma(r).$$
We conclude that $x_n\rightarrow r_\varepsilon$ entails $\sigma_n(x_n)\longrightarrow\sigma(r_\varepsilon)>0$, as required.%which proves (iii).

Finally, let $\gamma_\epsilon:=\inf_{|r-r_\varepsilon|<\epsilon}|\tilde{\varphi}'(r)|$, for $\epsilon>0$; since
$\tilde{\varphi}'(r_\varepsilon)<0$, for small $\epsilon$ we have $\gamma_\epsilon>0$. On the other hand, for every $n\geq 0$ it holds that
$$\tilde{\varphi}(r_{\varepsilon}^n)-\tilde{\varphi}_{t_n}(r_{\varepsilon}^n)=
\tilde{\varphi}(r_{\varepsilon}^n)-\tilde{\varphi}(r_\varepsilon)=-\tilde{\varphi}'(u_n)(r_{\varepsilon}^n-r_\varepsilon),$$
%\begin{eqnarray*}
%\tilde{\varphi}(r_{\varepsilon}^n)-\tilde{\varphi}_{t_n}(r_{\varepsilon}^n) & = & \tilde{\varphi}(r_{\varepsilon}^n)-
%\tilde{\varphi}(r_\varepsilon) \\
%& = & -\tilde{\varphi}'(u_n)(r_{\varepsilon}^n-r_\varepsilon),
%\end{eqnarray*}
for some $u_n\in(r_{\varepsilon},r_{\varepsilon}^n)$;
for the first equality we used the fact that $\tilde{\varphi}(r_\varepsilon)=\varepsilon=\tilde{\varphi}_{t_n}(r_{\varepsilon}^n)$,
while the second one follows by the mean value theorem. Consequently, for large $n$, satisfying $|r_{\varepsilon}^n-r_\varepsilon|<\epsilon$, we have
$$(r_{\varepsilon}^n-r_\varepsilon)\leq\gamma_\epsilon^{-1}\sup_{0\leq r\leq r_*}|\tilde{\varphi}_{t_n}(r)-\tilde{\varphi}(r)|=
o(n^{-\kappa});$$ this proves the claim and concludes the proof.
$\hfill\square$

\renewcommand{\tabcolsep}{3pt}
\begin{table}%[ht!]
\begin{tabular}{||c||c|c|c|c||c|c||}
\hline
$\log\varepsilon$       &   $-7$  &   $-8$  &   $-9$  &  $-10$  & $\bar{r}$ & $r_*$ \\
\hline\hline
$\hat{r}_\varepsilon$(A)& $0.619$ & $0.620$ & $0.621$ & $0.622$ &  $1.111$  & \:N/A\: \\
\hline\hline
$\hat{r}_\varepsilon$(B)& $0.502$ & $0.503$ & $0.504$ & $0.505$ &  $0.545$  & $0.5$ \\
\hline\hline
$\hat{r}_\varepsilon$(C)& $0.541$ & $0.542$ & $0.543$ & $0.544$ &  $0.545$  & \:N/A\: \\ \hline
\end{tabular}
\caption{\label{allR}\footnotesize{Critical threshold estimates ($\hat{r}_\varepsilon$) for various accuracy levels ($\varepsilon$) for three networks: A = Figure \ref{BD:fig}, B = Figure \ref{LK:fig} (Lu-Kumar), C = Figure \ref{LK:fig} (FCFS).}}
\end{table}

\begin{table}%[h!]
\begin{tabular}{||c||c|c|c|c|c|c||}
\hline
$\theta\setminus t$ &  $40$ &  $80$ & $100$ & $200$ &  $400$ & $1000$ \\
\hline\hline
\multicolumn{7}{|c|}{Network A}\\ \hline
%$0.00$	            & 1.000 & 1.000 & 1.000 & 1.000 & 1.000 & 1.000 \\ \hline
$0.11$	            & 0.865 & 0.865 & 0.865 & 0.865 & 0.865 & 0.865 \\ \hline
%$0.22$	            & 0.707 & 0.707 & 0.707 & 0.707 & 0.707 & 0.707 \\ \hline
$0.33$	            & 0.525 & 0.523 & 0.522 & 0.522 & 0.522 & 0.522 \\ \hline
%$0.44$	            & 0.334 & 0.315 & 0.313 & 0.309 & 0.308 & 0.308 \\ \hline
$0.55$	            & 0.174 & 0.135 & 0.125 & 0.104 & 0.092 & 0.086 \\ \hline
%$0.66$	            & 0.073 & 0.036 & 0.028 & 0.012 & 0.004 & 0.001 \\ \hline
$0.77$	            & 0.025 & 0.006 & 0.004 & 0.000 & 0.000 & 0.000 \\ \hline
%$0.88$	            &  0.03 &  0.01 &  0.00 & 0.00  & 1.00  &  1.00 \\ \hline
%$0.99$	            &  0.01 &  0.00 &  0.00 & 0.00  & 1.00  &  1.00 \\ \hline
%$1.11$	            &  0.01 &  0.00 &  0.00 & 0.00  & 1.00  &  1.00 \\ \hline
\hline\hline
\multicolumn{7}{|c|}{Network B}\\ \hline
$0.08$	            & 0.826 & 0.826 & 0.826 & 0.826 & 0.826 & 0.826 \\ \hline
%$0.16$	            & 0.642 & 0.641 & 0.641 & 0.641 & 0.641 & 0.641 \\ \hline
$0.24$	            & 0.455 & 0.451 & 0.451 & 0.450 & 0.450 & 0.450 \\ \hline
%$0.32$	            & 0.284 & 0.266 & 0.264 & 0.262 & 0.261 & 0.261 \\ \hline
$0.40$	            & 0.152 & 0.120 & 0.114 & 0.100 & 0.096 & 0.095 \\ \hline
%$0.48$	            & 0.070 & 0.039 & 0.033 & 0.018 & 0.010 & 0.006 \\ \hline
$0.56$	            & 0.028 & 0.009 & 0.006 & 0.001 & 0.000 & 0.000 \\ \hline
\hline\hline
\multicolumn{7}{|c|}{Network C}\\ \hline
$0.08$	            & 0.829 & 0.829 & 0.829 & 0.829 & 0.829 & 0.829 \\ \hline
%$0.16$	            & 0.655 & 0.655 & 0.655 & 0.655 & 0.655 & 0.655 \\ \hline
$0.24$	            & 0.484 & 0.483 & 0.483 & 0.482 & 0.482 & 0.482 \\ \hline
%$0.32$	            & 0.324 & 0.320 & 0.319 & 0.318 & 0.318 & 0.318 \\ \hline
$0.40$	            & 0.194 & 0.177 & 0.175 & 0.172 & 0.171 & 0.171 \\ \hline
%$0.48$	            & 0.102 & 0.078 & 0.072 & 0.062 & 0.058 & 0.057 \\ \hline
$0.56$	            & 0.047 & 0.026 & 0.021 & 0.012 & 0.007 & 0.003 \\ \hline
\end{tabular}
\caption{\footnotesize{Monotonicity shown in relation to the arrival rate($\theta$) and the time horizon ($t$) 
for three networks: A = Figure \ref{BD:fig}, B = Figure \ref{LK:fig} (Lu-Kumar), C = Figure \ref{LK:fig} (FCFS).}}
\label{allTheta}
%\caption{\label{LK-PRE}\footnotesize{Critical threshold estimates (top) and monotonicity (bottom) for the Lu-Kumar network,

%\caption{\label{LK-FIFO:tab}\footnotesize{Critical threshold estimates (top) and monotonicity (bottom) for the FCFS version ofthe Lu-Kumar network}.}
\vspace{-2mm}
\end{table}

\newpage

\bibliographystyle{IEEEtran}
\bibliography{BiB}

\renewcommand{\baselinestretch}{0.98}

\end{document}